\documentclass[reqno,centertags, 11pt]{amsart}
\usepackage{amsmath,amsthm,amscd,amssymb}
\usepackage{latexsym}
\newcommand{\bbE}{{\mathbb{E}}}
\newcommand{\bbR}{{\mathbb{R}}}
\newcommand{\bbD}{{\mathbb{D}}}
\newcommand{\bbT}{{\mathbb{T}}}
\newcommand{\bbZ}{{\mathbb{Z}}}

\newcommand{\calC}{{\mathcal{C}}}
\newcommand{\calF}{{\mathcal F}}
\newcommand{\calH}{{\mathcal H}}
\newcommand{\calL}{{\mathcal L}}
\newcommand{\calM}{{\mathcal M}}

\newcommand{\dott}{\,\cdot\,}

\newcommand{\lb}{\label}
\newcommand{\f}{\frac}

\newcommand{\ol}{\overline}
\newcommand{\ti}{\tilde  }
\newcommand{\wti}{\widetilde  }

\newcommand{\ac}{\text{\rm{ac}}}

\newcommand{\s}{\text{\rm{s}}}

\newcommand{\bi}{\bibitem}

\newcommand{\beq}{\begin{equation}}
\newcommand{\eeq}{\end{equation}}
\newcommand{\ba}{\begin{align}}
\newcommand{\ea}{\end{align}}

\newcounter{smalllist}

\newcommand{\bigtimes}{\mathop{\mathchoice%
{\smash{\vcenter{\hbox{\LARGE$\times$}}}\vphantom{\prod}}%
{\smash{\vcenter{\hbox{\Large$\times$}}}\vphantom{\prod}}%
{\times}%
{\times}%
}\displaylimits}

\DeclareMathOperator{\Real}{Re}

\allowdisplaybreaks
\numberwithin{equation}{section}

\newtheorem{theorem}{Theorem}[section]
\newtheorem*{t0}{Theorem}

\newtheorem*{p2.1}{Proposition 2.1}
\newtheorem{proposition}[theorem]{Proposition}

\theoremstyle{definition}

\theoremstyle{remark}

\newcommand{\abs}[1]{\lvert#1\rvert}

\begin{document}
\title[Aizenman's Theorem for OPUC]
{Aizenman's Theorem for Orthogonal Polynomials on the Unit Circle}
\author[B. Simon]{Barry Simon*}

\thanks{$^*$ Mathematics 253-37, California Institute of Technology, Pasadena, CA 91125.
E-mail: bsimon@caltech.edu. Supported in part by NSF grant DMS-0140592}

\date{September 27, 2004} 

\begin{abstract} For suitable classes of random Verblunsky coefficients, including 
independent, identically distributed, rotationally invariant ones, we prove that if 
\[
\bbE \biggl( \int\f{d\theta}{2\pi} \biggl|\biggl( \f{\calC + e^{i\theta}}{\calC-e^{i\theta}} 
\biggr)_{k\ell}\biggr|^p \biggr) \leq C_1 e^{-\kappa_1 \abs{k-\ell}} 
\]
for some $\kappa_1 >0$ and $p<1$, then for suitable $C_2$ and $\kappa_2 >0$, 
\[
\bbE \bigl( \sup_n \abs{(\calC^n)_{k\ell}}\bigr) \leq C_2 e^{-\kappa_2 \abs{k-\ell}} 
\]
Here $\calC$ is the CMV matrix. 
\end{abstract}

\maketitle

\section{Introduction} \lb{s1} 

This paper is a contribution to the theory of orthogonal polynomials on the 
unit circle (OPUC); for background on OPUC, see Szeg\H{o} \cite{Szb}, 
Geronimus \cite{GBk}, and Simon \cite{OPUC1,OPUC2}. Our goal here is to prove an 
analog of a result of Aizenman \cite{Aiz} for random discrete Schr\"odinger 
operators. Aizenman considers operators on $\ell^2 (\bbZ^\nu)$ of the form 
$h_\omega =h_0 + V_\omega$ where (\cite{Aiz} allows more general $h_0$ than 
this!) 
\[
(h_0 u)(n) = \sum_{\abs{j}=1} u(n+j) 
\]
and $V_\omega$ is the diagonal matrix whose matrix elements are independent 
identically distributed random variables. Aizenman's theorem states  

\begin{t0}[Aizenman \cite{Aiz}] \lb{AizT} Under suitable hypotheses on the distribution 
of $V$\!, if for some $[a,b]\subset\bbR$, 
\begin{equation} \lb{1.1} 
\int_a^b \bbE \bigl( \abs{[(h_\omega - E-i0)^{-1} P_{[a,b]}(h_\omega)]_{k\ell}}^p\bigr) 
\, dE \leq C_1 e^{-\kappa_1 \abs{k-\ell}} 
\end{equation} 
for some $0<p<1$ and $\kappa_1 >0$, then for some $\kappa_2>0$ and $C_2$, 
\begin{equation} \lb{1.2} 
\bbE \bigl( \sup_t \abs{[e^{-ith_\omega} P_{[a,b]} (h_\omega)]_{k\ell}}\bigr) 
\leq C_2 e^{-\kappa_2 \abs{k-\ell}} 
\end{equation} 
\end{t0} 

Aizenman's motivation was that Aizenman-Molchanov \cite{AM} had proven bounds of 
the form \eqref{1.1} (generally called Aizenman-Molchanov or fractional moment 
bounds) realizing that the key was to restrict $p$ to be less than $1$. From their 
bounds, they easily obtained spectral localization (i.e., pure point spectrum) by 
using the Simon-Wolff criterion \cite{SimWolff}. Aizenman was interested in 
\eqref{1.2} because it is a form of physical localization. It was used by 
del Rio et al.~\cite{Sim250} to obtain what is now the standard strong form of 
eigenfunction localization (SUDL) and by Minami \cite{Min} to prove Poisson 
distribution of the eigenvalues of $h_\omega$ restricted to a large box. Del Rio 
et al.~also simplified Aizenman's proof and slightly extended the result (so that  
the theorem we stated above is their form with some extra hypotheses dropped). 

To describe precisely the result we want to prove here, we need some preliminaries. 
Given a set of Verblunsky coefficients, $\{\alpha_j\}_{j=0}^\infty$ (see 
\cite[Section~1.5]{OPUC1}), one forms the CMV matrix, $\calC$ (see Cantero-Moral-Vel\'azquez 
\cite{CMV} or \cite[Chapter~4]{OPUC1}). We define, for $z\in\bbD$, the unit disk
\begin{equation} \lb{1.3} 
F_{k\ell}(z) = \biggl[ \f{\calC+z}{\calC-z}\biggr]_{k\ell} 
\end{equation} 
By Kolmogorov's theorem (see \cite{Kol} or Duren \cite[Section~4.2]{Duren}), $F_{k\ell}$  
lies in the Hardy spaces $H^p (\bbD)$ for $0<p<1$, so 
\begin{equation} \lb{1.4} 
F_{k\ell} (e^{i\theta}) \equiv \lim_{r\uparrow 1}\, F_{k\ell} (re^{i\theta}) 
\end{equation} 
exists for $\f{d\theta}{2\pi}$-a.e.~$\theta$ and has an integrable $p$-th power 
over $\partial\bbD$ for $p\in (0,1)$. 

Now let the $\alpha$'s be random variables which define a measure $d\Gamma$ on 
$\bigtimes_{j=0}^\infty \bbD\equiv \bbD^\infty$. For each $n=0,1,2,\dots$ and 
$\lambda\in\partial\bbD$, define $T_{n,\lambda} :\bbD^\infty \to \bbD^\infty$ by 
\begin{equation} \lb{1.4a}
\begin{aligned}
(T_{n,\lambda}(\alpha))_j &= \alpha_j && \qquad j=0,1,\dots, n-1 \\
&= \lambda\alpha_j && \qquad j = n, n+1, \dots 
\end{aligned}  
\end{equation}
and let $d\Gamma_{n,\lambda}(\alpha)=d\Gamma (T_{n,\lambda}(\alpha))$. We say 
$d\Gamma$ is strongly quasi-invariant if each $d\Gamma_{n,\lambda}$ is $d\Gamma$-absolutely 
continuous and $\sup_{n,\lambda} \|d\Gamma_{n,\lambda}/d\Gamma\|_\infty <\infty$. 
Clearly, if $d\Gamma$ is a product of rotation invariant measures (invariant i.i.d.'s), 
$d\Gamma$ is strongly quasi-invariant. We will discuss other examples in Section~\ref{s7}. 
Here is our main result: 

\begin{theorem}\lb{T1.1} Suppose $\{\alpha_j\}_{j=0}^\infty$ are random Verblunsky 
coefficients which are strongly quasi-invariant, and for some $p<1$ and $\kappa_1 >0$, 
\begin{equation} \lb{1.5} 
\bbE \biggl( \int_0^{2\pi} \, \abs{F_{k\ell} (e^{i\theta})}^p \, \f{d\theta}{2\pi}\biggr) 
\leq C_1 e^{-\kappa_1 \abs{k-\ell}}
\end{equation}
Then for suitable $\kappa_2 >0$ and $C_2$, 
\begin{equation} \lb{1.6} 
\bbE \big( \sup_n \, \abs{(\calC^n)_{k\ell}}\bigr) \leq C_2 e^{-\kappa_2 \abs{k-\ell}} 
\end{equation} 
\end{theorem} 

{\it Remarks.} 1. This result is interesting only because one can prove \eqref{1.5}. 
For certain cases of rotation invariant i.i.d.~$\alpha$'s, Stoiciu \cite{Stoi} has 
proven \eqref{1.5}. Indeed, I proved Theorem~\ref{T1.1} precisely to fill in a missing 
step in his program to prove Poisson distribution for the zeros of paraorthogonal 
polynomials with random Verblunsky coefficients. 

\smallskip
2. Kolmogorov's argument proves for any OPUC, any $k,\ell$, and $0<p<1$, 
\begin{equation} \lb{1.7} 
\int_0^{2\pi}\, \abs{F_{k\ell} (e^{i\theta})}^p \, \f{d\theta}{2\pi} 
\leq 2\cos \biggl( \f{p\pi}{2}\biggr) 
\end{equation} 
This means that if \eqref{1.5} holds for one $p\in (0,1)$, it holds for all such $p$ 
(with $\kappa_1$ dependent on $p$). 

\smallskip 
3. The $\sup_n$ is over $n=0, \pm 1, \pm 2, \dots$. 

\smallskip 
4. \eqref{1.6} is a strong statement about the structure of eigenfunctions of $\calC$ 
and of OPUC to be compared with the case $\alpha\equiv 0$ where $\sup_n 
\abs{(\calC^n)_{k\ell}} =1$. 

\smallskip
5. If \eqref{1.5} holds for $p=1$ (it cannot, as we will see in Remark~6!), \eqref{1.6} 
would be immediate since $2(\calC^n)_{k\ell}$ are the Taylor coefficients of $F_{k\ell}$. 
Thus, \eqref{1.5} $\Rightarrow$ \eqref{1.6} without recourse to the expectations. But, 
in general, for $p<1$, $H^p$ functions have Taylor coefficients that can grow as 
$o(n^{1/p-1})$ and no better (see Duren \cite[Chapter~6]{Duren}), so \eqref{1.5} 
$\Rightarrow$ \eqref{1.6} only holds because the $\sup \abs{(\calC^n)_{k\ell}}$ is 
averaged over a set of rank one perturbations. This is Aizenman's key discovery in 
\cite{Aiz}. 

\smallskip 
6. \eqref{1.5} cannot hold for $p=1$. Indeed, if $\sum_\ell (\int_0^{2\pi} 
\abs{F_{k\ell}(e^{i\theta})} \f{d\theta}{2\pi})^2 <\infty$ for fixed $k$, then 
$\sum_\ell \abs{\calC_{k\ell}^n}^2 \to 0$ by appealing to the dominated convergence 
theorem for sums and the Riemann-Lebesgue lemma which implies $(\calC^n)_{k\ell} 
\to 0$ if $\int \abs{F_{k\ell} (e^{i\theta})} \f{d\theta}{2\pi} <\infty$. But 
since $\calC^n$ is unitary, $\sum_\ell \abs{\calC_{k\ell}^n}^2 =1$. 

\smallskip
While part of our proof of Theorem~\ref{T1.1} follows the arguments in del Rio 
et al.~\cite{Sim250}, there are two novel aspects that prompted me to write 
this separate note. The first involves the theory of rank one perturbations. This 
theory is well-developed for selfadjoint operators (see \cite{AD,Scp38}), but 
I could not find any extensive theory for the unitary case when I wrote 
\cite{OPUC1,OPUC2}, so I developed the theory there (see \cite[Subsections~1.3.9, 
1.4.16, and Section~4.5]{OPUC1}). It turns out that a key formula needed here 
(see \eqref{2.12} below) is not in that presentation. 

Secondly, OPUC has a subtlety missing from the Schr\"odinger case. Namely, 
the relevant rank one perturbations of the Schr\"odinger operators are also 
Schr\"odinger operators, but the rank one perturbations of CMV matrices 
are not CMV matrices. Of course, as unitary matrices with a cyclic vector, 
these are unitarily equivalent to CMV matrices and, as we will see, a 
formula of Khrushchev \cite{Kh2000} even implies what the Verblunsky 
coefficients are for the new matrices. But we will need to know the form of 
the unitary, and this will require an illuminating calculation that should 
be useful in other contexts. 

In Section~\ref{s2}, we discuss some aspects of the theory of rank one  
perturbations of unitaries, which we use in Section~\ref{s3}, following 
\cite{Aiz,Sim250}, to compute explicit spectral representations. In Section~\ref{s4}, 
we use this to obtain a deterministic form of Aizenman's theorem that involves 
averaging under rank one perturbations, and in Section~\ref{s5}, we write these rank 
one perturbations in terms of CMV matrices. Section~\ref{s6} puts everything 
together to get Theorem~\ref{T1.1} and Section~\ref{s7} has some comments. 

\smallskip 
I would like to thank Mihai Stoiciu for useful discussions.  

\smallskip
\section{Rank One Perturbations of Unitaries Revisited} \lb{s2}  

Rank one perturbations of unitaries are best understood multiplicatively. 
Let $U$ be a unitary operator on a Hilbert space, $\calH$, and $\varphi 
\in\calH$ a unit vector. Let $P=\langle\varphi, \dott\rangle\varphi$ be the projection 
onto the multiples of $\varphi$. One defines for $\lambda\in\partial\bbD$, 
\begin{equation} \lb{2.1} 
U_\lambda = U(1-P) + \lambda U\!P = U[(1-P)+\lambda P] 
\end{equation} 
(if $\varphi$ is part of an orthonormal basis, $(1-P)+\lambda P$ is the diagonal 
matrix with $\lambda$ in $\varphi$ position and $1$ in all others). Thus 
\begin{align} 
U_\lambda \varphi & = \lambda U\varphi \lb{2.2}  \\
U_\lambda \psi &= U\psi \qquad\text{if } \psi\perp\varphi \lb{2.3} 
\end{align} 
which also defines $U_\lambda$. Note also that 
\begin{equation} \lb{2.3a} 
U_\lambda - U = (\lambda -1) U\!P 
\end{equation}

For $z\in\bbD$ and $\|\varphi\|=1$, define 
\begin{equation} \lb{2.4} 
F_\varphi (z) = \bigg\langle \varphi, \f{U+z}{U-z}\,\varphi\bigg\rangle 
\end{equation}
which is a Carath\'eodory function (i.e., $F(0)=1$ and $\Real F >0$ on $\bbD$). 
Since 
\begin{equation} \lb{2.5x} 
\bigg\langle \varphi, \f{U-z}{U-z}\, \varphi\bigg\rangle =1 
\end{equation}
we can solve for 
\begin{align} 
\langle\varphi, U(U-z)^{-1} \varphi\rangle &= \tfrac12\, [F_\varphi(z) + 1] \lb{2.5} \\
\langle\varphi, (U-z)^{-1} \varphi\rangle &= \tfrac{1}{2z}\, [F_\varphi(z) -1] \lb{2.6} 
\end{align} 

Notice \eqref{2.5} and \eqref{2.6} hold if $U$ is replaced by $U_\lambda$ and $F_\varphi$ 
by $F_\varphi^\lambda(z)$ given by \eqref{2.4} with $U$ replaced by $U_\lambda$. By 
the second resolvent equation, 
\begin{align} 
(z-U_\lambda)^{-1} U\varphi &= (z-U)^{-1} U\varphi + (z-U)^{-1} (U_\lambda -U) 
(z-U_\lambda)^{-1} U\varphi \notag \\
&= [1+(\lambda -1) (\varphi, (z-U_\lambda)^{-1} U\varphi)](z-U)^{-1} U\varphi \lb{2.7} 
\end{align}
on account of \eqref{2.3a}. 

Taking an inner product of \eqref{2.7} with $\varphi$ and using \eqref{2.5}/\eqref{2.6} for 
$U$ and $U_\lambda$ let us solve for $F_\varphi^\lambda (z)$ in terms of $F_\varphi (z)$. 
The result (see \cite[Subsection~1.4.16]{OPUC1}) is expressed most succinctly via the function 
$f$ defined by 
\begin{equation} \lb{2.8} 
F_\varphi(z) = \f{1+zf(z)}{1-zf(z)} 
\end{equation} 
for then the inner product of $\varphi$ with \eqref{2.7} implies 
\begin{equation} \lb{2.9} 
F_\varphi^\lambda (z) = \f{1+\lambda^{-1} zf(z)}{1-\lambda^{-1} zf(z)} 
\end{equation} 

Notice, by \eqref{2.5}/\eqref{2.6}, that 
\begin{align}
\langle\varphi, U_\lambda (U_\lambda -z)^{-1} \varphi\rangle &= \f{1}{1-\lambda^{-1} zf(z)} \lb{2.10} \\
\langle\varphi, (U_\lambda -z)^{-1} \varphi\rangle &= \f{\lambda^{-1} f(z)}{1-\lambda^{-1} zf(z)} \lb{2.11} 
\end{align} 

Thus far, the formulae are identical to what is in \cite{OPUC1}. What is new here is to note 
that \eqref{2.7} says that as a vector in $\calH$, for each $z\in\bbD$, $\lambda\in\partial\bbD$, 
$(U_\lambda -z)^{-1}U_\lambda \varphi$ is a multiple of $(U-z)^{-1} U\varphi$, so for any 
$\psi\in\calH$, 
\begin{equation} \lb{2.12} 
\f{\langle\psi, (U_\lambda-z)^{-1} U_\lambda\varphi\rangle}
{\langle\varphi, (U_\lambda-z)^{-1} U_\lambda\varphi\rangle} 
= \f{\langle\psi, (U-z)^{-1} U\varphi\rangle}{\langle\varphi, (U-z)^{-1} U\varphi\rangle} 
\end{equation}
In particular, by \eqref{2.10} and $\langle\psi, (U_\lambda-z)(U_\lambda-z)^{-1}
\varphi\rangle = \langle \psi,\varphi\rangle$, we see  

\begin{proposition}\lb{P2.1} If $\psi\perp\varphi$, then 
\begin{equation} \lb{2.13} 
\bigg\langle \psi, \f{U_\lambda +z}{U_\lambda -z}\,\varphi\bigg\rangle = 
\f{1-zf(z)}{1-\lambda^{-1} zf(z)}\, \bigg\langle \psi, \f{U+z}{U-z}\, \varphi\bigg\rangle 
\end{equation} 
\end{proposition} 

{\it Remark.} \eqref{2.12} is an analog of (3.2) of Aizenman \cite{Aiz}.  

\smallskip 
\section{The Spectral Representation} \lb{s3} 

We add two extra assumptions to our analysis of rank one perturbations of unitaries. 
First, we suppose $\varphi$ is cyclic for $U$\!, that is, $\{U^k\varphi\}_{k=-\infty}^\infty$ 
spans $\calH$, in which case it is easy to see that $\varphi$ is cyclic for $U_\lambda$. 
The spectral theorem then implies there are spectral measures, $d\mu_\lambda$, on $\partial 
\bbD$ defined by 
\begin{equation} \lb{3.1} 
\int \f{e^{i\theta}+z}{e^{i\theta}-z}\, d\mu_\lambda (\theta) = F_\varphi^\lambda (z) 
\end{equation}
for $z\in\bbD$ and unique unitary maps $\calF_\lambda : \calH \to L^2 (\partial\bbD, d\mu_\lambda)$, 
so that 
\begin{equation} \lb{3.2} 
(\calF_\lambda U_\lambda \psi)(z) =z (\calF_\lambda \psi)(z) \qquad 
\calF_\lambda\varphi \equiv 1
\end{equation} 
In particular, if $U_\lambda \eta_{z_0} =z_0 \eta_{z_0}$ for some $z_0\in\partial\bbD$ 
and $0\neq \eta\in\calH$, then, by \eqref{3.2}, 
\begin{equation} \lb{3.3} 
(\calF_\lambda \psi)(z_0) = \f{\langle \eta_{z_0},\psi\rangle}{\langle \eta_{z_0}, \varphi\rangle}  
\end{equation} 

Our second assumption, following \cite{AD,SimWolff}, is that for a.e.~$e^{i\theta_0}\in 
\partial\bbD$, 
\begin{equation} \lb{3.4} 
G(e^{i\theta_0}) \equiv \int \f{d\mu (\theta)}{\abs{e^{i\theta_0}-e^{i\theta}}^2} <\infty 
\end{equation} 
By arguments in \cite{SimWolff}, it is easy to see that if $\{\psi_j\}_{j=0}^\infty$ is a 
basis for $\calH$, then \eqref{3.4} at $\theta_0$ is equivalent to 
\begin{equation} \lb{3.5} 
\sum_{j=0}^\infty \, \lim_{r\uparrow 1}\, \biggl| \bigg\langle \psi_j, 
\biggl( \f{U+re^{i\theta_0}}{U-re^{i\theta_0}}\biggr) \varphi\bigg\rangle\biggr|^2 <\infty  
\end{equation} 
Moreover, if \eqref{3.4} holds, then 
\begin{equation} \lb{3.6} 
\lim_{r\uparrow 1}\, e^{i\theta_0} f(re^{i\theta_0}) =\lambda_0\in\partial\bbD 
\end{equation}
and $z_0$ is an eigenvalue of $U_{\lambda_0}$. Since (spectral averaging, due to 
Golinskii-Nevai \cite{GN} in this setting) 
\begin{equation} \lb{3.7} 
\int_0^{2\pi} (d\mu_{e^{i\varphi}}(\theta))\, \f{d\varphi}{2\pi} = \f{d\theta}{2\pi}  
\end{equation} 
\eqref{3.4} for a.e.~$e^{i\theta_0}$ implies $U_\lambda$ has pure point spectrum for 
a.e.~$\lambda$. (These facts are all explained in \cite[Sections~10.1 and 10.2]{OPUC2}.) 

\begin{proposition} \lb{P3.1} If $\varphi$ is cyclic and \eqref{3.4} holds at $\theta_0$, 
let $\lambda_0$ be given by \eqref{3.6}. Then for any $\psi\perp\varphi$, 
\begin{equation} \lb{3.8} 
(\calF_{\lambda_0}\psi)(z_0) = \ol{(1-\lambda_0)z_0} \, \, \, 
\ol{\langle\psi, (U-z_0)^{-1} \varphi\rangle}
\end{equation}
where $z_0 =e^{i\theta_0}$ and \eqref{3.8} is shorthand for $\lim_{r\uparrow 1} 
\ol{\langle\psi, (U-rz_0)^{-1}\varphi\rangle}$, which it is asserted exists. 
\end{proposition} 

\begin{proof} We use \eqref{2.12} for $z=rz_0$ and $\lambda =\lambda_0$. By the spectral 
theorem since $z_0$ is an eigenvalue of $U_{\lambda_0}$, 
\[ 
\lim_{r\uparrow 1}\, (1-r)(U_{\lambda_0} -rz_0)^{-1} U_{\lambda_0} \varphi 
\equiv \eta_{z_0} 
\]
is an eigenvector for $U_{\lambda_0}$ with eigenvalue $z_0$ or it is zero. Since \eqref{3.6} 
holds, we have that 
\begin{align*} 
\langle\varphi,\eta_{z_0}\rangle &= \lim_{r\uparrow 1}\, \f{1-r}{2} \, 
[F_\varphi^{\lambda_0} (rz_0)+1] \\
&= \lambda_0 \biggl[ \f{1-r}{\left.\lambda_0 - [zf(z)]\right|_{z=rz_0}}\biggr] 
\qquad \text{by \eqref{2.9}} \\
&= \f{\lambda_0 z_0}{(1-\lambda_0 z_0)^2}\, (-G(z_0))^{-1} 
\end{align*} 
by (10.1.7) of \cite{OPUC2}. This is nonzero by the assumption $G(z_0)<\infty$. One 
can also use cyclicity of $\varphi$ to conclude that $\eta_{z_0}\neq 0$. Thus 
$\eta_{z_0} \neq 0$ and so, by \eqref{3.3} and \eqref{2.12}, 
\begin{equation} \lb{3.9} 
(\calF_\lambda \psi)(z_0) = \lim_{r\uparrow 1}\, 
\f{\ol{\langle\psi, (U-rz_0)^{-1} U\varphi\rangle}}
{\,\ol{\langle\varphi, (U-rz_0)^{-1} U\varphi\rangle}\,} 
\end{equation} 
By \eqref{2.10} and $rz_0 f(rz_0)\to\lambda_0$, 
\begin{equation} \lb{3.10} 
\lim_{r\uparrow 1}\, \langle\varphi, (U-rz_0)^{-1} U\varphi\rangle = (1-\lambda_0)^{-1} 
\end{equation}
and since $(\psi, (U-z)^{-1} (U-z)\varphi)=0$, we have 
\begin{equation} \lb{3.11} 
\langle\psi, (U-rz_0)^{-1}U\varphi\rangle = rz_0 \langle\psi, (U-rz_0)^{-1}\varphi\rangle 
\end{equation} 
Thus, we see that \eqref{3.9}--\eqref{3.11} implies \eqref{3.8}. 
\end{proof} 

{\it Remark.} In fact, \eqref{3.8} holds whenever $d\mu$ is purely singular and for all 
$\lambda_0\neq 1$. This is because $d\mu$ purely singular implies $d\mu_\lambda$ is purely 
singular and \eqref{3.9} can be replaced by Poltoratskii's theorem \cite{Pol,JL}, 
which says that for any complex Borel measure $\eta$ on $\partial\bbD$ and any $g\in L^1 
(\partial\bbD, d\eta)$, we have, for almost any $e^{i\theta_0}$ with respect to $d\eta_\s$ 
(but not for $d\eta_\ac$), that 
\[
\lim_{r\uparrow 1}\, \f{[\int \f{e^{i\theta}+re^{i\theta_0}}{e^{i\theta}-re^{i\theta_0}} 
f(\theta)\, d\eta(\theta)]}{[\int \f{e^{i\theta}+re^{i\theta_0}}{e^{i\theta}-re^{i\theta_0}}\, 
d\eta(\theta)]} = f(\theta_0) 
\]

\smallskip
\section{Deterministic Form of Aizenman's Theorem} \lb{s4}

We now follow Aizenman \cite{Aiz} and del Rio et al.~\cite{Sim250}. Under the assumption 
that $\varphi$ is cyclic and $G(z_0)<\infty$ for a.e.~$z_0$ in $\partial\bbD$, we have for 
$\psi\perp\varphi$ that for a.e.~$\lambda_0$ that 
\begin{equation} \lb{4.1} 
\abs{(\calF_{\lambda_0}\psi)(z_0)} \leq \bigl| \lim_{r\uparrow 1}\, 
\langle \varphi, (U+rz_0)(U-rz_0)^{-1}\psi\rangle \bigr| 
\end{equation}
since $\psi\perp\varphi$ implies $\langle \varphi, (U+z)(U-z)^{-1}\psi\rangle = 
2z \langle \varphi, (U-z)^{-1}\psi\rangle$. \eqref{4.1} holds for all eigenvalues 
of $U_{\lambda_0}$ and so, for a.e.~$z_0$ w.r.t.~$d\mu_{\lambda_0}$ if $U_{\lambda_0}$ is 
pure point.  

Since $\calF_\lambda$ is a unitary operator, we have 
\begin{equation} \lb{4.2} 
\int \abs{(\calF_{\lambda_0}\psi)(z)}^2 \, d\mu_{\lambda_0}(z) = \|\psi\|^2 
\end{equation}
Moreover, since $\calF_{\lambda_0} U\calF_{\lambda_0}^{-1} =z$ and $\calF_{\lambda_0} 
\varphi\equiv 1$, 
\begin{align} 
\abs{\langle \varphi, U_{\lambda_0}^n\psi\rangle} 
&= \biggl| \int z^n \calF_{\lambda_0}\psi(z)\, d\mu_{\lambda_0}(z)\biggr| \notag \\ 
&\leq \int \abs{\calF_{\lambda_0}\psi(z_0)}\, d\mu_{\lambda_0}(z) \lb{4.3} 
\end{align} 

We conclude, using H\"older's inequality and \eqref{4.2}:  

\begin{proposition}\lb{P4.1} If $\varphi$ is cyclic and $G(z_0)<\infty$ for a.e.~$z_0$, then 
for a.e.~$\lambda_0$ and any $0<p<1$, we have 
\begin{equation} \lb{4.4} 
\sup_n\, \abs{\langle \varphi, U_{\lambda_0}^n \psi\rangle} \leq 
\biggl[ \int \lim_{r\uparrow 1}\, \abs{\langle\varphi, (U+rz)(U-rz)^{-1}\psi\rangle}^p\, 
d\mu_{\lambda_0}(z)\biggr]^{1/(2-p)} 
\end{equation} 
\end{proposition} 

\begin{proof} We have 
\[
\int \abs{g}\, d\mu_\lambda \leq \biggl( \int \abs{g}^2 \, d\mu_\lambda\biggr)^{(1-p)/(2-p)} 
\biggl( \int \abs{g}^p\, d\mu_\lambda\biggr)^{1/(2-p)} 
\]
since $(1-p)/(2-p) + 1/(2-p) =1$ and $2(1-p)/(2-p) + p/(2-p)=1$, and H\"older's 
inequality says $q\to \log (\int \abs{g}^q\, d\mu)^{1/q}$ is convex. \eqref{4.4} 
follows by taking $g=\calF_{\lambda_0} \psi$ and using \eqref{4.1}, \eqref{4.2}, and 
\eqref{4.3}. 
\end{proof} 

Since $2-p>1$, for any probability measure $d\nu$, $\int h^{1/(2-p)}\, d\nu \leq 
(\int h\, d\nu)^{1/(2-p)}$ by H\"older's inequality. Thus writing $\lambda_0 = 
e^{i\eta_0}$ and integrating \eqref{4.4} with $d\eta_0/2\pi$, we find, using 
\eqref{3.7}, that 

\begin{theorem}[Deterministic Aizenman's Theorem] \lb{T4.2} If $\varphi$ is cyclic, 
$\psi\perp\varphi$, and $G(z_0)<\infty$ for a.e.~$z_0\in\partial\bbD$, then for any $0<p<1$, 
\begin{equation} \lb{4.5} 
\int \f{d\eta}{2\pi}\, \sup_n\, \abs{\langle \varphi, U_{e^{i\eta}}^n \psi\rangle} 
\leq \biggl( \int \lim_{r\uparrow 1}\, \abs{\langle\varphi, (U+re^{i\theta})(U-re^{i\theta})^{-1} 
\psi\rangle}^p \, \f{d\theta}{2\pi} \biggr)^{1/(2-p)} 
\end{equation} 
\end{theorem}

\smallskip
\section{Rank One Perturbations of CMV Matrices} \lb{s5} 

We now specialize to $U=\calC (\{\alpha_j\}_{j=0}^\infty)$, a CMV matrix (see \cite{CMV} or 
\cite[Chapter~4]{OPUC1}), and $\varphi=\delta_n$, the vector with $1$ in the $n$-th position 
$(n=0,1,2,\dots$). We need a notation for diagonal matrices with diagonal matrix elements 
$\lambda$, $\lambda^{-1}$, or $1$. $D(\lambda^k (\lambda^{-1})^\ell (1\lambda)^\infty)$ 
will denote the diagonal matrix with $k$ $\lambda$'s, $\ell(\lambda^{-1})$'s, and then  
alternating $1$ and $\lambda$. We will also use the maps $T_{n,\lambda}$ of \eqref{1.4a}. 
Here is the main result: 

\begin{theorem}\lb{T5.1} Define $U_n$ for $n=0,1,2,\dots$ by 
\begin{align} 
U_{2k-1} &= D(1^{2k} (1\lambda)^\infty) \lb{5.1} \\
U_{2k} &= D(\lambda^{2k} (1\lambda)^\infty) \lb{5.2} 
\end{align} 
Then
\begin{equation} \lb{5.3} 
U_n\calC (T_{n,\lambda^{-1}} (\alpha)) U_n^{-1} = \calC(\alpha) \Delta_n (\lambda) 
\end{equation}
where 
\begin{equation} \lb{5.4} 
\Delta_n (\lambda) = D(1^n \lambda 1^\infty) 
\end{equation} 
\end{theorem} 

{\it Remarks.} 1. Since $\calC(\alpha) \Delta_n (\lambda)$ is unitary with $\delta_0$ cyclic, 
it is unitarily equivalent to some $\calC (\ti\alpha)$. Since $\calC(\alpha)\restriction 
\{\delta_j\}_{j=0}^{n-1} = \calC(\ti\alpha)\restriction \{\delta_j\}_{j=0}^{n-1}$, it is 
easy to see that $\ti\alpha_k =\alpha_k$ for $k=0,1,\dots, n-1$. By Khrushchev's formula 
(\cite[Theorem~9.2.4]{OPUC2}), the spectral measure for $\calC(\alpha)$ and vector 
$\delta_n$ has Schur function $\Phi_n (z;\alpha_0, \dots, \alpha_{n-1}) \Phi_n^* 
(z;\alpha_0, \dots, \alpha_{n-1})^{-1} f(z;\alpha_n, \alpha_{n+1}, \dots)$. By 
\eqref{2.9}, $\ti\calC(\alpha)$ thus has Schur function which is this times 
$\lambda^{-1}$, so by Khrushchev's formula again, $f(z; \ti\alpha_n, \ti\alpha_{n+1}, 
\dots) = \lambda^{-1} f(z;\alpha_n, \alpha_{n+1}, \dots) = f(z;\lambda^{-1} \alpha_n, 
\lambda^{-1} \alpha_{n+1}, \dots)$. We conclude $\ti\alpha_j = \lambda^{-1} \alpha_j$ 
for $j\geq n$. \eqref{5.3} goes beyond this by making the unitary equivalence explicit. 

\smallskip
2. The case $n=0$ is essentially Theorem~4.2.4 of \cite{OPUC1}. 

\begin{proof} First, some preliminaries. We use $\oplus$ for direct sum, normally of 
$2\times 2$ matrices but sometimes of a $1\times 1$ followed by $2\times 2$, in which 
case we write $\mathbf{1}_1$ or $\lambda\mathbf{1}_1$ so that, for example, if 
\begin{equation} \lb{5.3x} 
v(\lambda) = \begin{pmatrix} 1 & 0 \\ 0 & \lambda \end{pmatrix} \qquad 
\ti v (\lambda) = \begin{pmatrix} \lambda & 0 \\ 0 & 1 \end{pmatrix}
\end{equation} 
then 
\begin{align} 
U_{2k-1} &= \underbrace {\mathbf{1}_2 \oplus \mathbf{1}_2 \oplus \cdots \oplus\mathbf{1}_2}_{k \text{ times}} \,
\oplus \, v(\lambda) \oplus v(\lambda)\oplus \cdots \lb{5.4x}  \\
& =  \mathbf{1}_1 \oplus \underbrace{\mathbf{1}_2 \oplus\cdots\oplus \mathbf{1}_2}_{k\text{ times}}\, 
\oplus \, \ti v(\lambda) \oplus \ti v(\lambda) \oplus \cdots \lb{5.5}
\end{align}

If 
\[ 
\Theta(\alpha) = \left( \begin{array}{rr} \bar\alpha & \rho \\ \rho & -\alpha \end{array} \right) 
\]
with $\rho = (1-\abs{\alpha}^2)^{1/2}$, then (see \cite[Theorem~4.2.5]{OPUC1}) 
\begin{align}
\calC(\alpha) &= \calL(\alpha) \calM(\alpha) \lb{5.6} \\ 
\calL(\alpha) &= \Theta(\alpha_0) \oplus \Theta(\alpha_2) \oplus \Theta(\alpha_4) \oplus \cdots \lb{5.7} \\
\calM(\alpha) &= \mathbf{1}_1 \oplus \Theta(\alpha_1)\oplus \Theta(\alpha_2) \oplus \cdots \lb{5.8} 
\end{align} 

Notice next that (note $v(\lambda)$ not $v(\lambda)^{-1}$ in both places!) 
\begin{align} 
v(\lambda)\Theta(\lambda^{-1} \alpha) v(\lambda) &= \lambda\Theta(\alpha) \lb{5.9} \\
\ti v(\lambda)^{-1} \Theta(\lambda^{-1}\alpha) \ti v(\lambda)^{-1} &= \lambda^{-1} \Theta(\alpha) \lb{5.10} 
\end{align} 

We now turn to the proof of \eqref{5.3} for $n=2k-1$. By \eqref{5.4x}, \eqref{5.7}, and \eqref{5.9}, 
\begin{equation} \lb{5.11} 
U_{2k-1} \calL (T_{2k-1, \lambda^{-1}}(\alpha)) U_{2k-1} = \lambda\calL(\alpha) W
\end{equation}
where 
\begin{equation} \lb{5.12} 
W=D((\lambda^{-1})^{2k} 1^\infty) 
\end{equation} 
Now 
\begin{align} 
WU_{2k-1}^{-1} &= D((\lambda^{-1})^{2k} (1\lambda^{-1})^\infty) \notag \\
&= \lambda^{-1} \mathbf{1}_1 \oplus 
\underbrace{ \lambda^{-1} \mathbf{1}_2 \oplus\cdots\oplus\lambda^{-1}\mathbf{1}_2}_{(k-1) \text{ times}} \, 
\oplus \,\ti v(\lambda)^{-1} \oplus \ti v(\lambda)^{-1} \oplus \cdots  \lb{5.13} 
\end{align} 
since $\lambda^{-1} (1\lambda^{-1})^\infty = (\lambda^{-1} 1)^\infty$. Thus, by \eqref{5.10} 
and \eqref{5.8}, 
\begin{equation} \lb{5.14} 
WU_{2k-1}^{-1} \calM(T_{2k-1, \lambda^{-1}}(\alpha)) U_{2k-1}^{-1} W = 
\lambda^{-1} \calM(\alpha) D((\lambda^{-1})^{2k-1} 1^\infty) 
\end{equation}
We conclude 
\begin{align*} 
U_{2k-1} & \calC (T_{2k-1, \lambda^{-1}}(\alpha)) U_{2k-1}^{-1} \\
&= U_{2k-1} \calL(T_{2k-1,\lambda^{-1}}(\alpha))U_{2k-1} U_{2k-1}^{-1} 
\calM (T_{2k-1, \lambda^{-1}}(\alpha)) U_{2k-1}^{-1} \\
&= \lambda\calL(\alpha) WU_{2k-1}^{-1} \calM(T_{2k-1,\lambda^{-1}}(\alpha)) U_{2k-1}^{-1} 
WW^{-1} \quad\text{(by \eqref{5.11})} \\
&= \calL(\alpha)\calM(\alpha) D((\lambda^{-1})^{2k-1} 1^\infty) D((\lambda)^{2k}1^\infty) 
\quad\text{(by \eqref{5.14})} \\
&=\calC(\alpha) \Delta_n(\lambda) 
\end{align*} 
since $\Delta_n(\lambda) = (1^n \lambda 1^\infty)$. This proves \eqref{5.3} for $n=2k-1$. 

Now suppose $n=2k$. Then, by \eqref{5.9}, 
\begin{equation} \lb{5.15} 
U_{2k} \calL (T_{2k,\lambda^{-1}}(\alpha)) U_{2k} = \lambda\calL (\alpha )\wti W 
\end{equation} 
where 
\begin{equation} \lb{5.16} 
\wti W= D((\lambda)^{2k} 1^\infty) 
\end{equation}
Thus 
\begin{align} 
\wti W U_{2k}^{-1} &= D(1^{2k} (1\lambda^{-1})^\infty) \notag \\
&= \mathbf{1}_1 \oplus \underbrace{\mathbf{1}_2 \oplus\cdots\oplus\mathbf{1}_2}_{k \text{ times}}\,  
\oplus \, \ti v(\lambda)^{-1} \oplus \ti v(\lambda)^{-1} \oplus \cdots \lb{5.17} 
\end{align} 
This is the reason odd and even $n$ differ. \eqref{5.13} has $k-1$ $\mathbf{1}_2$'s and 
\eqref{5.17} has $k$ of them. By \eqref{5.10} and \eqref{5.8},  
\begin{equation} \lb{5.18} 
\wti WU_{2k}^{-1} \calM (T_{2k,\lambda^{-1}}(\alpha)) U_{2k}^{-1} \wti W = 
\lambda^{-1} \calM(\alpha) D(\lambda^{2k-1} 1^\infty)  
\end{equation}
It follows that 
\begin{align*} 
U_{2k} \calC (T_{2k, \lambda^{-1}}(\alpha))U_{2k} 
&= U_{2k} \calL (T_{2k,\lambda^{-1}}(\alpha)) U_{2k} U_{2k}^{-1} \calM (T_{2k,\lambda^{-1}}(\alpha)) 
U_{2k}^{-1} \\
&= \lambda \calL(\alpha) \wti WU_{2k}^{-1} \calM (T_{2k,\lambda^{-1}}(\alpha)) U_{2k}^{-1} 
\wti W \wti W^{-1} \\
&= \calL (\alpha) \calM(\alpha) D(\lambda^{2k+1} 1^\infty) D((\lambda^{-1})^{2k} 1^\infty) \\
&= \calC(\alpha) \Delta_n (\lambda)
\qedhere
\end{align*} 
\end{proof} 

\smallskip
\section{Proof of Theorem~\ref{T1.1}} \lb{s6} 

We are now ready to put it all together: 

\begin{proof}[Proof of Theorem~\ref{T1.1}] By \eqref{3.5} and the exponential decay 
in \eqref{1.5}, $G(z_0)<\infty$ for a.e.~$\alpha$, so Theorem~\ref{T4.2} applies for 
a.e.~$\alpha$. By \eqref{4.5} with $\varphi =\delta_k$ and $\psi =\delta_m$ and 
$U=\calC(\alpha)$, 
\begin{equation} \lb{6.1} 
\int \f{d\varphi}{2\pi}\, \sup_n\, \abs{[(\calC(\alpha)\Delta_k (e^{i\varphi}))^n]_{km}} 
\leq \biggl( \int \abs{F_{km} (e^{i\theta})}^p \, \f{d\theta}{2\pi}\biggr)^{1/2} 
\end{equation} 
Now use Theorem~\ref{T5.1} and the fact that $U$ is diagonal and unitary to see  
\[
\abs{[(\calC(\alpha) \Delta_k (\lambda))^n]_{km}} = 
\abs{[\calC(T_{k,\lambda^{-1}}(\alpha))^n]_{km}} 
\]
so that 
\begin{equation} \lb{6.2} 
\int \f{d\varphi}{2\pi} \, \sup_n\, \abs{[\calC(T_{k,e^{-i\varphi}}(\alpha))^n]_{km}} 
\leq \biggl( \int \abs{F_{km}(e^{i\theta})}^p \f{d\theta}{2\pi} \biggr)^{1/2}   
\end{equation}
Take expectations of both sides. Use the quasi-invariance to write 
\[
\bbE\bigl( \sup_n \abs{[\calC (T_{k,e^{-i\varphi}}(\alpha))^n]_{k\ell}}\bigr) 
\geq C^{-1} \bbE \bigl( \sup_n \abs{[C(\alpha)^n]_{k\ell}}\bigr) 
\]
for a constant $C$ independent of $k,\ell,\varphi$. Use H\"older's inequality to bring 
$\bbE$ inside $(\cdot)^{1/(2-p)}$. The result is 
\[
\bbE \bigl( \sup_n \abs{(\calC^n)_{k\ell}}\bigr) \leq C \biggl( \int 
\bbE (\abs{F_{k\ell} (e^{i\theta})}^p)\, \f{d\theta}{2\pi}\biggr)^{1/2} 
\]
which shows \eqref{1.5} implies \eqref{1.6}. 
\end{proof} 

\smallskip
\section{Remarks} \lb{s7} 

Some closing remarks: 

\smallskip 
1. It is not hard to prove a local version of this theorem where $\int_0^{2\pi}$ in 
\eqref{1.5} is replaced by $\int_a^b$ and an extra $P_{(a,b)}(\calC)$ is added in 
\eqref{1.6}. This might be useful in the quasi-invariant case, but in the i.i.d.~rotation 
invariant case, $\bbE (\abs{F_{k\ell} (e^{i\theta})}^p)$ is $\theta$-invariant and 
so the integral has exponential decay for $(a,b)$ if and only if it does for $(0,2\pi)$. 

\smallskip
2. If $d\rho$ is a rotation quasi-invariant measure on $\bbD$ and $d\gamma$ one on 
$\partial\bbD$, and if $\alpha_0, \bar\alpha_0\alpha_1, \bar\alpha_1 \alpha_2, \dots$ 
are independent random variables with $\alpha_0$ $d\rho$-distributed and each 
$\bar\alpha_j\alpha_{j+1}$, $d\rho$-distributed, then this measure is 
quasi-invariant. It would be interesting to do localization theory (both spectral 
and dynamic) for this model. 

\smallskip
3. It would be interesting to know if \eqref{1.6} implies \eqref{1.5}.

\bigskip

\end{document}